\documentclass{article}
\usepackage{amsmath}
\usepackage[normalem]{ulem}
\usepackage{amssymb}
\usepackage{amsthm}
\usepackage{xcolor}
\usepackage[utf8]{inputenc}
\usepackage[letterpaper, margin=1in]{geometry}

\usepackage{booktabs}
\usepackage{adjustbox}
\usepackage{multirow}
\usepackage{soul}

\usepackage{algorithm}
\usepackage{algpseudocode}
\usepackage[
pdftitle={QCF},
pdfcreator={pdflatex},
pdfsubject={preprint},
hyperindex = {true},
colorlinks = {true},
linkcolor = {blue},
citecolor = {blue}
]{hyperref}

\usepackage{float}
\usepackage[font=scriptsize,skip=5pt]{caption}
\usepackage[subrefformat=parens]{subcaption}

\DeclareMathOperator*{\argmin}{arg\,min}

\newcommand{\xhdr}[1]{\vspace{5pt} \noindent {\textbf{#1}}}
\makeatletter
\DeclareRobustCommand\onedot{\futurelet\@let@token\@onedot}
\def\@onedot{\ifx\@let@token.\else.\null\fi\xspace}

\usepackage{bm}
\usepackage{bbm}
\newcommand{\xb}{\bm{x}}

\newcommand{\sig}{\sigma}
\newcommand{\pl}{\partial}
\newcommand{\trace}{{\rm trace}}

\newcommand{\footremember}[2]{%
    \footnote{#2}
    \newcounter{#1}
    \setcounter{#1}{\value{footnote}}%
}
\newcommand{\footrecall}[1]{%
    \footnotemark[\value{#1}]%
} 

\title{
Quantum-inspired variational algorithms for partial differential equations: Application to financial derivative pricing
}

\author{%
  Tianchen Zhao\footremember{UM}{Department of Mathematics, University of Michigan, Ann Arbor, MI 48109}%
  \and Chuhao Sun\footrecall{UM}
  \and Asaf Cohen\footrecall{UM}%
  \and James Stokes\footremember{FI}{Flatiron Institute, Simons Foundation, New York, NY 10010}%
  \and Shravan Veerapaneni\footrecall{UM}\, \footrecall{FI}%
  }

\date{}
\begin{document}
\maketitle

\begin{abstract}
    Variational quantum Monte Carlo (VMC) combined with neural-network quantum states offers a novel angle of attack on the curse-of-dimensionality encountered in a particular class of partial differential equations (PDEs); namely, the real- and imaginary time-dependent Schr\"odinger equation. In this paper, we present a simple generalization of VMC applicable to arbitrary time-dependent PDEs, showcasing the technique in the multi-asset Black-Scholes PDE for pricing European options contingent on many correlated underlying assets.
\end{abstract}
\tableofcontents
\section{Introduction}

The field of stochastic variational algorithms has undergone dramatic recent developments based on two relatively recent, albeit independent, developments: ({\em i}) the availability of near-term quantum computers \cite{preskill2018quantum} and ({\em ii}) the existence of scalable stochastic algorithms for training deep neural networks. The former development has motivated a new research direction called variational quantum algorithms (VQAs) \cite{cerezo2021variational}, in which the stochastic variational character of the algorithms render them suitable for noisy intermediate-scale quantum computers \cite{preskill2018quantum}. The latter holds promise to accelerate a host of scientific computing problems including general-purpose solvers of partial differential equations (PDEs) using physics-informed neural-networks (PINNs) \cite{raissi2019physics} and has already made substantial headway in solving the time-(in)dependent Schr\"{o}dinger equation in high dimensions using variational quantum Monte Carlo (VMC) with neural-network quantum states \cite{carleo2017solving}. It is noteworthy that there exist many shared parallels between the fields of VQAs, VMC \cite{stokes2020quantum, stokes2022numerical} and PINNs. In particular, both VMC and VQAs hinge on the concept of adaptive stochastic estimation of the time-independent and time-dependent variational principles originally due to Rayleigh-Ritz and McLachlan \cite{mclachlan1964variational}, respectively. VQAs and VMC differ essentially in the choice of parametrized quantum state and the associated adaptive sampling strategy used to estimate quantum expectation values. PINNs can be viewed as an alternative approach to McLachlan's variational principle. They gained significant traction outside the quantum physics literature and have been advocated for solving general PDEs that appear in other areas of science and engineering. This approach 
involves solving a non-convex optimization problem for the PDE residual in order to determine the space-time development of the state variable. A large number of variants of this proposal have been subsequently put forward targeting PDEs either in many spatial dimensions (e.g., \cite{sirignano2018dgm, wang2021spacetime}) or parametric PDEs in low dimensions (e.g., \cite{li2020fourier}). Deep backward stochastic differential equation methods \cite{weinan2017, Han2018, beck2019, hure2020} are another class of deep learning approaches that have been applied to parabolic PDEs that arise in mathematical finance.

In this paper, we introduce a generalization of McLachlan's variational principle applicable to a wide variety of time-dependent PDEs and propose a variational quantum Monte Carlo (VMC) stochastic approximate solution method utilizing autoregressive neural-network quantum states. A closely related algorithm has been recently introduced in the VQA literature \cite{fontanela2021quantum, alghassi2021variational}, which was motivated by the problem of solving linear PDEs using digital quantum computers.  Ref.~\cite{alghassi2021variational} argued for an approximate solution concept, with respect to which exponential quantum speedup could be achieved for state evolution, which however is overwhelmed by auxiliary costs of state preparation and extraction of properties of the state (see \cite[Section III.C]{alghassi2021variational} for details). These complications of read-out and state preparation are not relevant to the VMC-based solver, however, since the VMC computing model permits efficient queries to arbitrary probability amplitudes. Interestingly, our generalization of McLachlan's variational principle is closely related to the neural Galerkin method recently put forward in \cite{bruna2022neural}, which also introduced a different neural-network-based stochastic solution method.  One small difference compared to \cite{bruna2022neural} is that we choose to work on a predefined mesh, which provides a simple method to implement non-trivial boundary conditions necessary for financial applications. The use of a mesh is not mandatory, however, and this paper paves the way toward a mesh-free generalization. Indeed, in the final stages of preparation of this article, Ref.~\cite{reh2022variational} appeared, which proposes a mesh-free flow-based solution of probabilistic PDEs such as the Fokker-Planck equation. Unlike \cite{reh2022variational}, however, our approach does not require that the state variable of the underlying PDE corresponds to a probability density.

The speedup obtainable by the approach advocated here has practical applications in overcoming the curse-of-dimensionality in high-dimensional PDEs, particularly in situations where fine-grained information about the state variable is required, such as gradients with respect to the independent variables. One such situation is quantum many-body physics, where kinetic energy observables depend on second-order spatial derivatives. Another example is the pricing and hedging of contingent claims in multi-asset financial markets.
The approach is applicable to general time-dependent PDEs expressible in first-order form, although we only consider the inhomogeneous linear case in the numerical experiments. Specifically, we chose to focus on the multi-asset Black-Scholes PDE for pricing European contingent claims because of its well-known relation with the  time-dependent Schr\"{o}dinger equation (TDSE) as well as its importance in computational finance. This work opens the door to free-boundary problems necessary for pricing
American contingent claims.

The organization of the paper is as follows: In section \ref{sec:background} we generalize the McLachlan variational principle to general time-dependent PDEs in a model-agnostic manner, which is applicable in the purely classical or quantum setting. In section \ref{sec:modeling} we describe the modeling assumptions involved in the passage from a time-dependent PDE to an neural quantum state-based solution of McLachlan's variational principle.
The remainder of the paper is dedicated to numerical experiments, focusing on the problem of financial derivative pricing, which suffers from a curse-of-dimensionality. Section \ref{sec:experiments} in particular provides numerical confirmation in the case of correlated diffusions, and section \ref{sec:option} describes the application of these results to option pricing in the Black-Scholes pricing framework.

\section{Theory}\label{sec:background}
\subsection{Generalities of McLachlan's variational principle}
McLachlan's variational principle \cite{mclachlan1964variational} is an example of a time-dependent variational principle (TDVP) that approximates the solution of the TDSE by evolution within a space of parametrized trial functions. TDVPs for the TDSE have been devised for tensor network states \cite{haegeman2011time}, neural-network quantum states (NQS) \cite{carleo2017solving} and parametrized quantum circuits \cite{stokes2020quantum, barison2021efficient}. Ref.~\cite{stokes2022numerical} studied TDVPs through the lens of information geometry providing a unified perspective applicable to variational quantum algorithms (VQAs) and variational quantum Monte Carlo with normalized neural-network quantum states. In the following section, we further generalize TDVPs to include general time-dependent PDEs and establish additional connections with the VQA and numerical analysis literature. 
Given a time-dependent PDE for the state variable $u(t,x) \in \mathbb{C}$ with $(t,x) \in [0,T] \times \Omega$, together with a choice of parametrized functions of the spatial variable $\{ u_\theta : \Omega \longrightarrow \mathbb{C} \, | \, \theta \in \mathbb{R}^p \}$, the output of a TDVP is parametrized curve $\gamma : [0,T] \longrightarrow \mathbb{R}^p$ in the space of variational parameters such that $u_{\gamma(t)}$ optimally describes $u(t,\cdot)$ in some distance metric $\operatorname{dist}(\cdot,\cdot)$ for all $t \in [0,T]$. Consider the initial value problem for a general time-dependent PDE of the form,
\begin{align}
    \partial_t u(t,x) & = \mathcal{F}(t,x,u) \\
    u(0,x) & = u_0 \in L^2(\Omega;\mathbb{C}) \enspace ,
\end{align}
with prescribed boundary conditions on $\partial \Omega$. In order to simplify exposition, in this section we avoid complications of boundary conditions by assuming either $\Omega = \mathbb{R}^d$ with suitable decay at infinity or $\Omega = \mathbb{T}^d$. In practice we replace the spatial domain by a mesh $\widehat{\Omega} \subset \Omega$, thereby approximating square-integrable functions by square summable vectors. The imposition of boundary conditions on the spatial mesh is then achieved via the use of source functions. Let $\Phi^{t}_s$ denote the time evolution map for state variable such that $\Phi^t_u\circ\Phi^u_s = \Phi^t_s$ and in particular $u(t,\cdot)=\Phi^t_0(u_0)$. 
Given an initial parameter vector $\theta_0 \in \mathbb{R}^p$ and a step size $\delta t > 0$, define a sequence of parameter vectors $( \theta_k )_{k \in \mathbb{N} }$ by the following iteration
\begin{align}
    \theta_{k+1}
    & := \argmin_{\theta \in \mathbb{R}^p} 
    \left[
    \operatorname{dist}\left(
    \Phi^{(k+1)\delta t}_{k \delta t}(u_{\theta_k})
    ,
    u_{\theta}\right)
    \right]. \label{e:projected}
\end{align}
In quantum physics applications, a suitable distance metric is the Fubini-Study metric, as previously argued both in VMC \cite{carleo2017solving, stokes2022numerical} and in VQA \cite{stokes2020quantum, barison2021efficient} literature. In this work we choose $\operatorname{dist}(\cdot,\cdot)$ to be the Euclidean norm. Rather than solving the discrete-time dynamical system \eqref{e:projected} directly, we consider the limit of infinitesimal step size $\delta t \longrightarrow 0$ in which it reduces to the following system of ordinary differential equations (ODEs)\footnote{See appendix \ref{app:evolution} for a derivation. Similiar ODEs have been introduced in the neural Galerkin method \cite{bruna2022neural}.} with initial condition $\gamma(0)=\theta_0$:
\begin{equation}\label{e:general}
    M(\gamma(t)) \, \gamma'(t) = V(t,\gamma(t))
\end{equation}
where
\begin{equation}
    M_{ij}(\theta) 
    := \operatorname{Re}\left[\left\langle \frac{\partial u_\theta}{\partial \theta_i} \middle| \frac{\partial u_\theta}{\partial \theta_j} \right\rangle\right] \enspace , \quad \quad 
    V_i(t,\theta) := \operatorname{Re}\left[\left\langle \frac{\partial u_\theta}{\partial \theta_i} \middle| \mathcal{F}(t,u_\theta) \right\rangle\right]
\end{equation}
and where $\langle \cdot | \cdot \rangle $, $\Vert \cdot \Vert_2$ denote the standard inner product and the induced norm for $L^2(\Omega,\mathbb{C})$, respectively. Although the matrix $M_{ij}$ is necessarily positive semi-definite, it may be degenerate, reflecting the possibility of multiple minima in the optimization problem \eqref{e:projected}. Thus, regularization techniques are generally required in order to obtain a well-posed system of ODEs.

\subsection{Analogy with finite element approximations}
Before discussing our autoregressive neural-network quantum state implementation, it is useful to orient within scientific computing literature by showing that the formalism 
shares close parallels with
the classic finite element method applied to the linear inhomogeneous case $\mathcal{F}(t,x,u) = \mathcal{L}u(t,x) + f(t,x)$. Given a set of real basis functions $\{\varphi_i\}_{i=1}^m$, define the variational family consisting of a weighted superposition,
\begin{equation}
    u_\theta(x) = \sum_{i=1}^m \theta_i \varphi_i(x) \enspace ,
\end{equation}
where $\theta \in \mathbb{R}^m$ is assumed here.
Substituting the above into \eqref{e:general} one finds the following ODE determining the dynamics of the weights,
%
\begin{equation}\label{e:fem}
    M \gamma'(t) = - K \gamma (t) + f(t) \enspace , \quad \quad M_{ij} := \langle \varphi_i | \varphi_j \rangle \enspace , \quad \quad K_{ij} := - \langle \varphi_i | \mathcal{L}u_j \rangle \enspace , \quad \quad f_i(t) := \langle \varphi_i | f(t) \rangle \enspace ,
\end{equation}
which can be recognized as the standard discrete linear systems arising in finite element methods, with $M$ and $K$ the mass- and stiffness-matrix, respectively. 
Recall that since finite element techniques use non-overlapping elements to discretize the spatial domain and employ basis functions with localized support, the problem size scales as $m \sim  p^d$ in $d$ spatial dimensions, where $p$ is the average number of elements per dimension. 
Thus, the curse-of-dimensionality arises from need to perform increasingly high-dimensional, albeit sparse, linear algebra in order to solve the linear system \eqref{e:fem}. The approach advocated in the following section, in contrast, overcomes the curse-of-dimensionality by representing the solution vector in terms of an autoregressive neural-network quantum state.

\subsection{Neural-network quantum state implementation}
In order to overcome the curse of dimensionality with respect to the spatial dimension $d$, which is inherent in the utilization of a $d$-dimensional mesh $\widehat{\Omega} \subset \Omega \subseteq \mathbb{R}^d$, we utilize stochastic estimation combined with autoregressive assumptions. In particular, we take inspiration from both the variational quantum algorithm introduced in \cite{alghassi2021variational} as well variational quantum Monte Carlo using autoregressive NQS \cite{sharir2020deep, hibat2020recurrent}, by parametrizing the solution of the PDE as $u_\theta = \alpha \, \psi_{\beta}$ where $\psi_{\beta}$ is a unit-normalized NQS with variational parameters $\beta$ and $\alpha > 0$ is a scale factor, whose time-dependence must be determined from the evolution equations along with $\beta$. The variational parameters thus consist of the augmented parameter vector $\theta := (\log\alpha,\beta) \in \mathbb{R}^{p+1}$ where $\beta \in \mathbb{R}^p$ is an unconstrained vector representing the weights and biases of the neural network. Plugging the rescaled ansatz into the evolution equations \eqref{e:general}, one obtains an augmented system of first-order, non-linear ordinary differential equations determining the time-dependence of the augmented vector $\theta$. The overlaps defining $M$ and $V$ are estimated using the VMC importance sampling technique, where the probability density is chosen to be the modulus-squared wavefunction $|\psi_{\beta}(x)|^2$. 



\section{Numerical implementation}\label{sec:modeling}
In this section, we provide detailed modeling assumptions required to implement the algorithm. After describing a mesh-based encoding of the state variable into a multi-qubit state, we then introduce the model architecture and computation mechanisms including the forward pass and the autoregressive sampling process. Pre-training is necessary in order to satisfy the initial condition of the PDE, and we adopt the standard approach~\cite{sirignano2018dgm} by updating the model iteratively on batches of randomly selected mesh points.
Finally, we describe the stochastic estimation procedure used to evolve the variational trial state using McLachlan's variational principle.

\subsection{Conversion to meshed form}
Rather than working in the continuum, we assume spatial discretization of the function $u(t,\cdot)$ on a regular grid $\widehat{\Omega}$ embedded in the $d$-dimensional domain $\Omega  = [a_1, b_1] \text{\texttimes} ... \text{\texttimes} [a_d, b_d]$ with a total of $2^n$ grid points. Without loss of generality, we assume $\Omega $ is a regular hypercube satisfying $|b_1-a_1|=...=|b_d-a_d|$, and that the mesh size along each axis is $2^{n/d}$. In order to represent the state of the discretized field in terms of the state of an $n$-qubit system, we assign each computational basis state $|k_1,\ldots,k_n \rangle \in \mathbb{C}^{2^n}$ with $(k_1,\ldots,k_n) \in \{0,1\}^n$ to a linear index defined by $k=\sum_{i=1}^n k_i 2^i$  and then unravel the linear index to indices along each of $d$ axes yielding a $d$-tuple $(\bar k_1,\ldots,\bar k_d) \in \{0, \ldots, n/d-1\}^d$ defined by
\begin{equation}
    \sum_{i=1}^n \bar k_i (2^{\frac{n}{d}})^{i-1}=k.
\end{equation}
The mesh point $\xb^k \in \widehat\Omega$ corresponding to the linear index $k \in \{0,1,...,2^n-1\}$ is thus
\begin{equation}
    \xb^k=(a_1+\bar k_1 \Delta x,...,a_d+\bar k_d \Delta x)
\end{equation}
where $\Delta x=|b_1-a_1| / 2^{\frac{n}{d}}$. Henceforth, we do not distinguish the index $k$, the binary string $(k_1,\ldots,k_n )$ and the corresponding coordinate $\xb^k \in \widehat\Omega$. The digitized representation of the state of the PDE at time $t$ is thus the following unnormalized $n$-qubit state
\begin{align}\label{e:|u(t)>}
    | u(t) \rangle = \sum_{k\in\widehat{\Omega}}u(\xb^k, t) | k_1,\ldots,k_n \rangle \enspace 
\end{align}
and likewise the parametrized approximation $|u_\theta\rangle$ of $|u(t)\rangle$ is given by,
\begin{equation}\label{e:utheta}
| u_\theta \rangle
=
\sum_{k\in\widehat{\Omega}}u_\theta(\bm{x}^k) | k_1,\ldots,k_n \rangle
\end{equation}
where $u_\theta(\cdot)$ is a function defined on $\widehat{\Omega} \subset \Omega$.
\subsection{Autoregressive assumption and sampling}
In order to obtain an expressive family of trial functions $u_\theta : \widehat{\Omega} \longrightarrow \mathbb{C}$ which furthermore admits an efficient stochastic estimation procedure, we express $u_\theta$ as multiple of a unit-normalized neural-network quantum state $\psi_\beta : \widehat{\Omega}\longrightarrow \mathbb{C}$ with variational parameters $\beta \in \mathbb{R}^p$,
\begin{equation}
    u_{\theta}(\bm{x}^k) = 
    \alpha \, \psi_{\beta}(k_1,\ldots,k_n) \enspace , \quad \quad \Vert \psi_\beta \Vert_{\widehat{\Omega}} = 1 \enspace .
\end{equation}
A simple method to ensure unit-normalization and efficient sampling is to follow the work of MADE~\cite{germain2015made} which exploited a masked version of a fully connected layer, where some connections in the computational path are removed in order to satisfy the autoregressive properties, In particular, if we assume a choice of variables such that the state variable is strictly positive\footnote{This can be considered as a special case of the complex-valued case \cite{sharir2020deep,hibat2020recurrent}.} $u_\theta(x) > 0$, then a suitable choice of unit-normalized function is
\begin{align}
    \psi_{\beta}(k_1,\ldots,k_n) = \prod_{i=1}^n \sqrt{p_{\beta,i}(k_i|k_{i-1},\ldots,k_1)} \enspace ,
\end{align}
where $p_{\beta,i}(\cdot|k_{i-1},\ldots,k_1)$ is the parametrized conditional probability distribution for the $i$th bit.

\subsection{Pre-training}
Before we perform the time evolution, initial variational parameters $\theta_0$ must be selected in order to match the variational function $u_{\theta_0}(x)$ with the choice of initial condition $u_0(x)$ for $x \in \widehat{\Omega}$ in the spatial domain.
This can be achieved via optimization of the following objective function using stochastic gradient descent,
\begin{align}
     J(\alpha_0,\beta_0) = \big\Vert \alpha_0 | \psi_{\beta_0}\rangle - | u_0 \rangle \big\Vert_{\widehat{\Omega}}^2 \enspace .
\end{align}
\subsection{Evolution}
In this section we discuss the details of the stochastic estimation of $M$ and $V$ necessary to evolve the augmented parameter vector $\theta = (\log\alpha,\beta) \in \mathbb{R}^{p+1}$. For simplicity we only consider the affine case
\begin{equation}\label{e:affine}
\mathcal{F}(t,x,u) = \mathcal{L}u(t,x) + f(t,x) \enspace .
\end{equation}
Consider the decomposition,

\begin{align}
    M &= 
    \left[\begin{array}{c|c } 
    	M_{00} & M_{0,1:p} \\ 
    	\hline 
    	M_{1:p,0} & M_{1:p,1:p} 
    \end{array}\right], \quad \quad  V = \left[\begin{array}{c} 
    	V_{0} \\ 
    	\hline 
    	V_{1:p}
    \end{array}\right].\nonumber \\
\end{align}
It is convenient to introduce the following helper functions. In particular, define the Born probability distribution
$\rho_\beta(x) \in [0,\infty)$, the wavefunction score $\sigma_\beta(x) \in \mathbb{C}^n$ and the local energy $l_\theta(t,x) \in \mathbb{C}$ as follows,
\begin{equation}\label{e:local_energy}
    \rho_\beta(x) := |\psi_\beta(x)|^2 \enspace , \quad \quad \sigma_\beta(x) := \frac{\nabla_\beta \psi_\beta(x)}{\psi_\beta(x)} \enspace , \quad\quad l_\theta(t,x) := \frac{(\mathcal{L}\psi_\beta)(x)}{\psi_\beta(x)} + \frac{f(t,x)}{\alpha} \enspace .
\end{equation}
By straightforward calculus, we obtain
\begin{equation}
    M_{00} = 
    \alpha^2 \enspace , \quad \quad
    M_{1:p,0} = M_{0,1:p} =
    \alpha^2 \operatorname{Re}\left[\underset{x \sim \rho_\beta}{\mathbb{E}} \sigma_\beta(x)\right], \quad \quad
    M_{1:p,1:p} = 
    \alpha^2 \operatorname{Re}\left[\underset{x \sim \rho_\beta}{\mathbb{E}} \overline{\sigma_\beta(x)} \, \sigma_\beta(x)^T\right], \\
\end{equation}
and
\begin{equation}
    V_{0} = 
    \alpha^2\operatorname{Re}\left[\underset{x \sim \rho_\beta}{\mathbb{E}} l_\theta(t,x)\right] \enspace , \quad\quad
    V_{1:p} =
    \alpha^2\operatorname{Re}\left[\underset{x \sim \rho_\beta}{\mathbb{E}} \overline{\sigma_\beta(x)}\, l_\theta(t,x)\right].
\end{equation}
The expectation values over $x$ are approximated using Monte Carlo sampling.
In practice, the batch of randomly generated samples is represented in the form of a buffer $\mathcal{B}=\{x_i\}_{i=1}^B$ that stores the unique samples in the batch and a counter $\mathcal{C}=\{c_i\}_{i=1}^B$ that records the number of occurrences of each of the corresponding samples. Expectation values are then approximated by sums of the following form,
\begin{align}
    \underset{x \sim \rho_\beta}{\mathbb{E}}[g(x)] \approx \frac{1}{\sum_{i=1}^B c_i} \sum_{i=1}^B c_i g(x_i).
\end{align}

\xhdr{Computation of per-sample gradient.}
The computation for both $M$ and $V$ requires direct access to the per-sample gradients $\{\nabla_{\beta} \psi_{\beta}(x_i)\}_{i=1}^B$, which is typically not directly accessible during a traditional backward pass by some deep learning software. 
In order to avoid inefficient forward and backward passes for each of the $B$ samples, we utilize \texttt{BackPack} which collects the quantities necessary to compute the individual gradients and reuses them to compute the per-sample gradient without significant computational overhead.

\xhdr{Parallel extraction of matrix elements.}
Given a sample $x \in \widehat{\Omega}$, corresponding to a particular row-index of the operator $\widehat{\mathcal{L}} \in \mathbb{C}^{2^n \times 2^n}$, the calculation of the local energy $l_\theta(t,x)$ defined in \eqref{e:local_energy} involves determining the nonzero entries of $\widehat{\mathcal{L}}$ in that row. Fortunately, the structure of the PDE problem ensures that the location and values of these entries can be determined in $O(\text{poly}(n))$ time.
We exploit CPU parallelization to determine the nonzero row entries for each sample in the batch. Note that the maximal number of nonzero entries per row is determined in advance; e.g. $2n^2+1$ for diffusion operator with Dirichlet boundary conditions.

\xhdr{Boundary Conditions.} Boundary conditions are implemented on the grid using an appropriate choice of source function. In the case of Dirichlet conditions, for example, we choose the source function as follows
\begin{equation}
    f(t,x) = u(t,x) \mathbbm{1}_{\partial \Omega}(x)
\end{equation}
where $\mathbbm{1}_{\partial \Omega}$ denotes the indicator function for the set $\partial \Omega$.

\xhdr{Parameter update.}
The network can be trained by incrementing the parameters using a Euler scheme in the direction $\delta \theta \in \mathbb{R}^{p+1}$ given by the solution of the following linear system
\begin{align}
    M\delta \theta = V \delta t,
\end{align}
where $M,V$ are computed as discussed before. In practice, $M$ is usually ill-conditioned due to the fact that it's essentially a sum of $B$ rank one matrices. To stabilize the inverse operation, we consider the singular value decomposition $M=U \Sigma W^T$, and remove the diagonal values of $\Sigma$ smaller than a small threshold $\epsilon=10^{-12}$ to obtain $\Sigma_r$ of size $r \text{\texttimes} r$, as well as $U_r, W_r$ of size $(p+1) \text{\texttimes} r$. The direction vector is therefore approximated by $W_r\Sigma_r^{-1}U_r^TV\, \delta t$.


\section{Diffusion with Gaussian Initializations}\label{sec:experiments}

In this section, we show various numerical experiments  demonstrating the convergence and run time analysis of our proposed approach. 
Here, we consider the $d$-dimensional heat equation for which $\mathcal{L}=D \,\nabla\cdot\nabla$  in \eqref{e:affine} (diffusion constant set   to $D=0.1$) with either periodic or Dirichlet boundary conditions on $\partial \Omega$. 
For benchmarking purposes, we employ a finite difference method with the standard central difference scheme to discretize $\mathcal{L}$ (formulas given in Appendix \ref{app:centraldiff}) and forward Euler for time-stepping. 

\subsection{Experimental setup}\label{subsec:set-ups}
For initialization we chose a discrete isotropic Gaussian, expressed in terms of the modified Bessel function $I_x(t)$ of integer order $x$,
\begin{equation}\label{e:init}
    u_0(x) = \prod_{i=1}^d e^{-t} I_{x_i}(t) \enspace ,
\end{equation}
where $t>0$ is a parameter controlling the width of the Gaussian, and $x$ ranges from 0 to $2^{n-1}$.
Pre-training was performed using Adam optimizer~\cite{kingma-iclr15} for 50k iterations with batch size $128$, $\beta_1=0.9$, $\beta_2=0.999$, and $\epsilon=10^{-8}$.
The learning rate was warm-started for the first $1/10$ total training iterations, then decayed by a factor of 10 at $3/7, 5/7$ of the total training iterations.

After completion of pre-training, the state was evolved for a total evolution time $T=1$ using a step size of $\delta t = 5\text{\texttimes}10^{-5}$ and a batch size of $B=1024$ for Monte Carlo estimation on each iteration. The time development of the state was compared to the result of Euler time-stepping the initial condition \eqref{e:init} using the same step-size. Due to the exponential scaling of the matrix $\widehat{\mathcal{L}}$, we only establish this baseline for the number of qubits $n\leq16$. All reported numerical results are averaged across exactly 5 seeds, each trained for 20k steps.





Table~\ref{tbl:benchmark_performance} reports the relative error between our method and forward Euler method over 20k iterations, for diffusion problem with Dirichlet boundary conditions. The relative error of the obtained solution is computed by comparing the norm of the space-time history of the variational state with the result of Euler development,
\begin{align}
     \textrm{error} := \frac{1}{T}\sum_{t \in \widehat{\mathcal{T}}}\frac{\big\Vert |u(t)\rangle - |u_{\gamma(t)}\rangle \big\Vert_{\widehat{\Omega}}}{\big\Vert |u(t)\rangle\big\Vert_{\widehat{\Omega}}} \enspace .
\end{align}
In general, the overall solution obtained using our method matches well with that obtained using forward Euler method. Their discrepancy increases as the dimensionality of the problem increases; nonetheless, within a tolerable threshold (e.g., less than $10\%$). One remedy is to simply increase the batch size, for example,  we show in Section~\ref{sec:batch_size} that the relative error improves from $15\%$ to around $3\%$ with a batch size that is ten times larger.

\begin{table*}[t]
\scriptsize
\centering
\caption{Average relative error of our method in comparison with forward Euler method over 20k iterations. The proposed algorithm solution is compared with a Euler forward method for the diffusion equation over 20k iterations. For forward Euler method, the time step is $5\text{\texttimes}10^{\text{\textminus}5}$ with a total time of $1$. The relative error is computed as $\frac{1}{T}\sum_{t=1}^T\frac{\lVert u(t,x) - f(x;\theta_t) \rVert}{\lVert u(t,x)\rVert}$.}
\label{tbl:benchmark_performance}
\begin{tabular}{lcccccc}
\toprule
\multirow{2}{*}{Operator} & \multirow{2}{*}{Boundary Condition} & \multirow{2}{*}{$n/d$} & \multicolumn{4}{c}{\# of Dimensions $d$} \\
\cmidrule(lr){4-7}
& & & 1 & 2 & 3 & 4 \\
\midrule
Diffusion & Dirichlet & 4 & $5.13 \text{\texttimes} 10^{\text{\textminus}3}$ & $7.92 \text{\texttimes} 10^{\text{\textminus}3}$ & $3.12 \text{\texttimes} 10^{\text{\textminus}2}$ & $1.47 \text{\texttimes} 10^{\text{\textminus}1}$ \\
Diffusion & Dirichlet & 5 & $2.91 \text{\texttimes} 10^{\text{\textminus}3}$ & $9.91 \text{\texttimes} 10^{\text{\textminus}3}$ & $7.24 \text{\texttimes} 10^{\text{\textminus}2}$ & - \\
\toprule
\end{tabular}
\end{table*}

\begin{table*}[t]
\scriptsize
\centering
\caption{Average running time over 2k iterations. The batch size used here is 500. Forward Euler method suffers from the exponential complexity, whereas our method, despite having an overhead running time, enjoys a polynomial scaling. Note that we cannot apply Euler for higher dimensions due to the memory constraint.}
\label{tbl:benchmark_running_time}
\begin{tabular}{lcccccccccccc}
\toprule
\multirow{2}{*}{Operator} & \multirow{2}{*}{Method} & \multirow{2}{*}{$n/d$} & \multicolumn{9}{c}{\# of Dimensions $d$} \\
\cmidrule(lr){4-12}
& & & 1 & 2 & 3 & 4 & 5 & 6 & 7 & 8 & 9 \\
\midrule
Diffusion & Euler & 4 & 0.024 & 0.137 & 1.559 & 305.92 & - & - & - & - & - \\
Diffusion & Euler & 5 & 0.031 & 0.225 & 68.827 & - & - & - & - & - & - \\
\midrule
Diffusion & Ours & 4 & 29.52 & 55.76 & 138.88 & 230.02 & 360.80 & 518.00 & 792.63 & 1214.09 & 1812.94 \\
Diffusion & Ours & 5 & 32.53 & 88.75 & 173.43 & 311.82 & 507.93 & 847.85 & 1355.16 & 2379.86 & 4123.61 \\
\toprule
\end{tabular}
\end{table*}

\subsection{Running Time Analysis}
Since we discretize the domain $\Omega$ using a grid of size $|\widehat{\Omega}|=2^n$, the time complexity of forward Euler method, expressed in terms of $n$ is $O(T \text{\texttimes} 2^{2n})$, where $T$ is the number of iterations. The proposed VMC algorithm, in contrast, scales as $O(T B\,\text{poly}(n))$, where $B$ is the batch size. In more detail, the forward pass scales as $O(T B n^2)$, and the sampling is $O(T B n^3)$, due to the sequential nature of the auto-regressive sampling process. This polynomial scaling comes at a price of the approximate nature of the time evolution step, and the implicit access to entries of the state vector compared to the forward Euler method which offers $O(1)$ lookup to entries of the state.

In Table~\ref{tbl:benchmark_running_time}, we report the average running time of both forward Euler method and our method for 2k iterations (one-tenth of the total running time). The batch size $B$ of our method is fixed to be 500. Note that we cannot apply Euler for higher dimensions due to the memory constraint. Although forward Euler method is effective for small-scale problems, its complexity suffers from the exponential growth with respect to dimension $d$. The computational cost of our method originates from four sources: sampling, forward pass, per-sample gradient computation (backward pass), and extraction of matrix element information.
All these sources contribute to the overhead time that causes our method to run slower in comparison with respect to forward Euler method for problems in lower dimensions. However, as the dimensionality increases, the run time of our method grows only at a polynomial rate. It can be seen from Table \ref{tbl:benchmark_running_time} that our method is already faster than forward Euler for problem sizes characterized by $n=16$ qubits.


\subsection{Convergence Visualization}

\begin{figure}[t]
\centering
\includegraphics[width=0.80\textwidth]{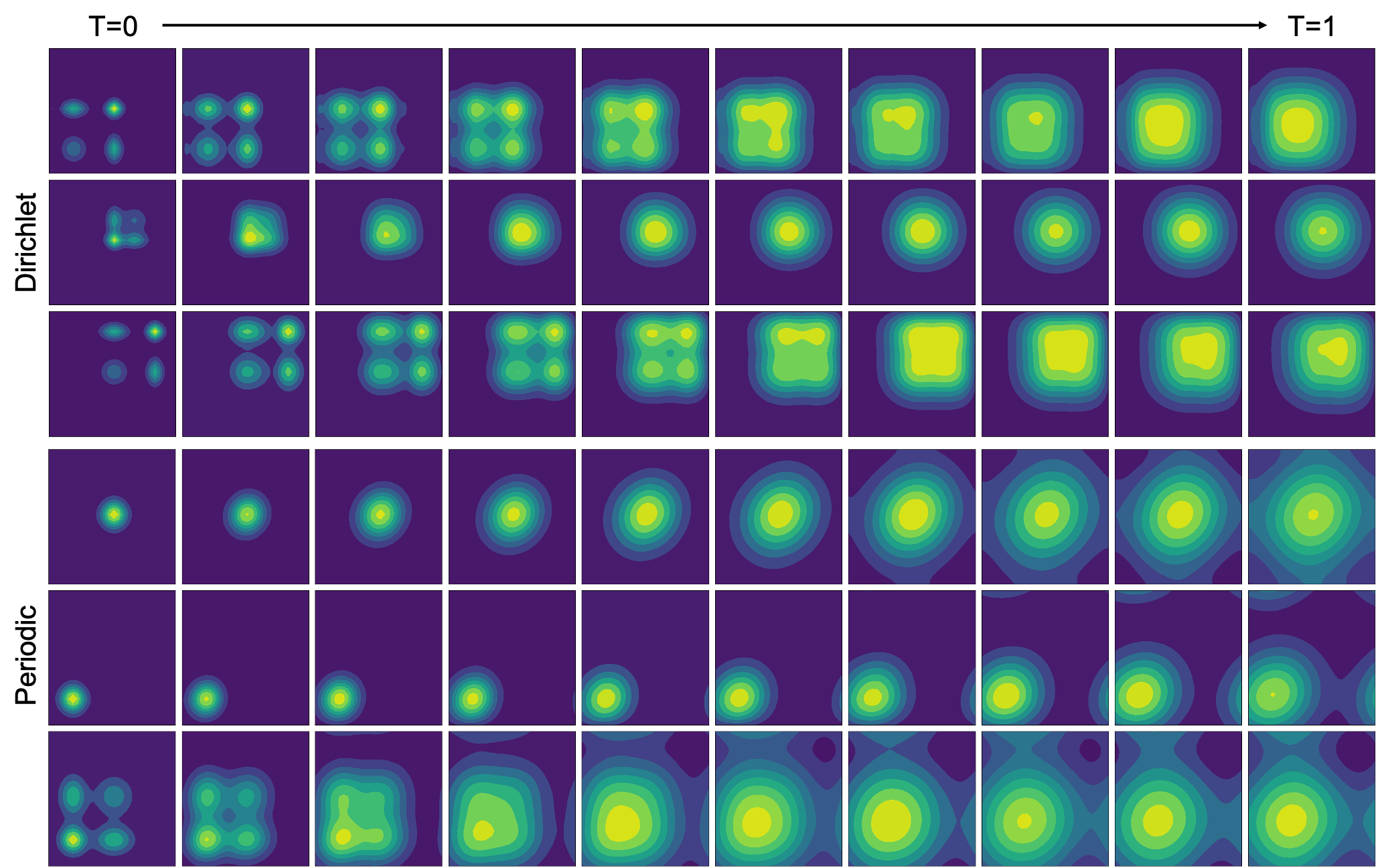}
\caption{Snapshots of the evolution obtained using our algorithm for diffusion equation with Dirichlet and periodic boundary conditions and different choices of initialization.}
\label{fig:snapshots}
\end{figure}

We provide snapshots of our method for 2D diffusion problems with periodic boundary conditions over 20k iterations. We run our algorithm with five distinct initializations and record the snapshot every 2k iterations.
In particular we do not employ any regularization techniques to enforce that the solution satisfies the boundary condition. It can be observed in Figure~\ref{fig:snapshots} that our method successfully  obeys the periodic boundary conditions.

\subsection{Ablation Study on Batch Size}
\label{sec:batch_size}
\begin{figure}[t]
\centering
\includegraphics[width=0.95\textwidth]{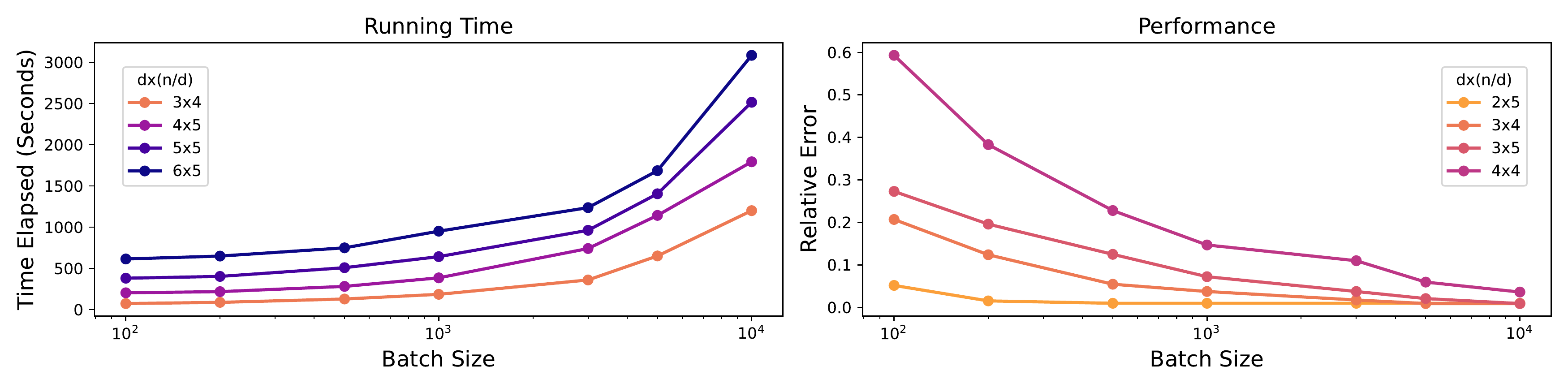}
\vspace{-2mm}
\caption{We report the running time and the average relative error with the forward Euler method over various batch sizes, for qubit sizes $n=d\times\frac{n}{d}$, where $d$ is the dimensionality of the problem. The running time grows with respect to the batch size and the problem size. On the other hand, the performance is greatly improved when training with larger batch sizes.}
\label{fig:batch_size}
\end{figure}

Recall that $M$ and $V$ are approximated with Monte Carlo sampling using batches of unique samples. Intuitively, a larger batch size yields a better approximation to the exact expectation value, thereby providing more accurate model updates. In this section, we study the effect of batch size on the performance of our method.
In the LHS of Figure~\ref{fig:batch_size}, the running time of our algorithm increases for both larger problem sizes and batch sizes. Note that the actual running time does not grow linearly with respect to the batch size in the plot due to the cache and parallelization.
In the RHS of Figure~\ref{fig:batch_size}, we report the average relative error between the forward Euler method and the VMC method with various batch sizes. Given a fixed problem size (e.g., 4 dimensions with 4 qubits per dimension), we observe a performance improvement by increasing the batch size, which verifies our hypothesis that increasing the batch size does effectively improve the performance. Given a fixed batch size, our method performs worse as the dimensionality of the problem increases. This result implies that we need a larger batch size to guarantee good performance for problems in higher dimensions.


\section{Option Pricing}\label{sec:option}
In this section, we explain how our results can be used to price options. We start with the well-known Feynman-Kac representation theorem and then describe how the numerical results obtained in Section \ref{sec:experiments} can be applied to the multi-dimensional Black-Scholes model. 

\subsection{The Feynmam-Kac representation theorem and option pricing}
The Feynman-Kac stochastic
representation formula links between a parabolic PDE and stochastic differential equations. Namely, consider a PDE of the form:
\begin{align*}
    \begin{cases}\frac{\partial u}{\partial t}(t,x)+\mu(t,x)\cdot\nabla_xu(t,x)+\frac12\trace[\sigma^TH\sigma](t,x)-r(t,x)u(t,x)+f(t,x)=0,\qquad (t,x)\in[0,T]\times\mathbb{R}^d,\\
    u(T,x)=\Phi(x),
    \end{cases}
\end{align*}
where $\mu:[0,T]\times\mathbb{R}^d\to\mathbb{R}^d$, $\sigma:[0,T]\times\mathbb{R}^d\to\mathbb{R}^{d\times d}$, $r,f:[0,T]\times\mathbb{R}^d\to\mathbb{R}$, and $\Phi:\mathbb{R}^d\to\mathbb{R}$ are given measurable functions and $H$ stands for the Hessian matrix of $u$:
\begin{align*}
   H_{ij}=\frac{\partial^2 u}{\partial x_i\partial x_j}.
\end{align*}
Then, under a second order integrability condition, the solution to this equation takes the form:
\begin{align}\label{eq:FKrep}
    u(t,x)=\mathbb{E}\Big[\int_t^Te^{-\int_t^\tau r(s,X_s)ds}f(\tau,X_\tau)d\tau+e^{-\int_t^Tr(s,X_s)ds}\Phi(X_T)\Big|X_t=x \Big],
\end{align}
where $(X_t)_{t\in[0,T]}$ is the solution to the following stochastic differential equation with $(W_t)_{t\in[0,T]}$ being a $d$-dimensional standard Wiener process:
\begin{align*}
    \begin{cases}
    dX_t=\mu(t,X_t)dt+\sigma(t,X_t)dW_t,\\
    X_t=x.
    \end{cases}
\end{align*}

One of the many applications of this representation is in derivative pricing. Consider multi-asset with price dynamics $(S_t)_{t\in[0,T]}$ given by the multi-dimension geometric form:
\begin{align*}
    \begin{cases}
    dS_t=\text{Diag}[S_t]\mu(t,S_t)dt+\text{Diag}[S_t]\sigma(t,S_t)dW_t,\\
    S_0=s_0>0,
    \end{cases}
\end{align*}
where $\text{Diag}[x]$ is a diagonal matrix with the vector $x$ on the diagonal, $W$ is $d$-dimensional standard Wiener process, 
$\mu:[0,T]\times[0,\infty)\to\mathbb{R}^d$ and $\sigma:[0,T]\times[0,\infty)\to\mathbb{R}^{d\times d}$ are measurable and satisfy basic conditions that lead to a unique strong solution to the stochastic differential equation above. Consider also the interest rate $r:[0,T]\times [0,\infty)\to\mathbb{R}$. 
A {\em risk-neutral measure}
is a probability measure $\mathbb{Q}$, such that the discounted price of the asset $(e^{-\int_0^tr(s,X_s)ds}X_t)_{t\in[0,T]}$ is a martingale under $\mathbb{Q}$. The dynamics of the multi-asset under $\mathbb{Q}$ are given by:
\begin{align*}
    \begin{cases}
    dS_t=r(t,S_t)S_t
    dt+\text{Diag}[S_t]\sigma(t,S_t)dW^\mathbb{Q}_t,\\
    S_0=s_0>0,
    \end{cases}
\end{align*}
where $W^{\mathbb{Q}}_\cdot:=W_\cdot-\int_0^\cdot (\mu(t,S_t)-r(t,S_t))dt$ is a standard Wiener process under $\mathbb{Q}$. 
Then, the price of a simple contingent claim with terminal payment $\Phi(S_T)$ at any given time $t\in[0,T]$, when the multi-asset prices are $S_t=s$ is given by $u(t,s)$, which satisfies the PDE:
\begin{align*}
    \begin{cases}\frac{\partial u}{\partial t}(t,s)+r(t,s)s\cdot\nabla_xu(t,s)+\frac12\trace[\sigma^TH\sigma](t,s)-r(t,s)u(t,s)=0,\qquad (t,s)\in[0,T]\times(0,\infty)^d,\\
    u(T,s)=\Phi(s),
    \end{cases}
\end{align*}
and can be expressed as the conditional $\mathbb{Q}$-expectation as follows:
\begin{align*}
    u(t,s)=\mathbb{E}^{\mathbb{Q}}\Big[e^{-\int_t^Tr(u,S_u)du}\Phi(S_T)\Big|S_t=s \Big].
\end{align*}

The most celebrated example is the Black-Scholes model, where the dynamics of the underlying asset follows a geometric Brownian motion and the interest rate is fixed. In dimension 1, the associated PDE admits an explicit solution, which is known as the {\it Black-Scholes formula}. Aside from the fact that the solution of the PDE is the contingent claim price, its partial derivative $\partial u/\partial x$ is used in order to construct a {\it delta-hedging} portfolio. Furthermore, its partial derivatives, known as the {\it Greeks} are used to construct a robust portfolio against small changes, when moving from the continuous- to the discrete-time world. 
\subsection{Option pricing in Black-Scholes model}
We now demonstrate how our algorithm applied to the multidimensional heat equation can be used to price options. 
For this we consider the Black-Scholes model that consists of a risk-free asset with a constant risk-free return $r>0$ and $d$ risky assets whose dynamics are given by
\begin{align*}
    dS^i_t = \mu_i S^i_t dt + \sig_i S^i_tdW_t^i,\qquad i=1,\ldots,d.
\end{align*}
The parameters $\mu_i$ and $\sig_i$, $i = 1,\dots, d$ are constants, and $\{W^i\}_{i=1}^d$ are Wiener processes with quadratic covariation $[W^i_t,W^j_t] = \rho_{ij}t$. 

Consider further a European option, whose payment at (the predetermined) expiry time $T>0$ is $\Psi(S_T^1,\dots,S^d_T)$, for some measurable function $\Psi$. Let $V$ be the conditional price of this option, i.e., $V(t,x_1,\dots,x_d)$ is the price for the option at time $t$, given that $S_t^i = x_i$ for $i = 1,\dots,d$. It is well-known that $V$ satisfies the following Black-Scholes PDE:
\begin{align*}
    \begin{cases}
    \frac{\pl V}{\pl t} + \sum_{i=1}^drx_i\frac{\pl V}{\pl x_i} + \frac12\sum_{i=1}^d \sig^2_ix_i^2\frac{\pl^2 V}{\pl x_i^2} + \sum_{i\ne j} \frac 12 \rho_{ij}\sig_i\sig_jx_ix_j\frac{\pl^2 V}{\pl x_i\pl x_j} - rV = 0,\qquad (t,x)\in[0,T]\times(0,\infty)^d,\\
    V(T,x_1,\ldots,x_d) = \Psi(x_1,\dots,x_d).
    \end{cases}
\end{align*}

Following Guillaume \cite{gui2019}, we may reduce this $n$-dimensional equation to the $n$-dimensional standard heat equation. To this end, set
\begin{align*}
    u(t,y_1,\ldots,y_n)=V(T-t,e^{\sigma_1y_1},\ldots,e^{\sigma_n y_n})e^{\sum_{i=1}^n-a_i\sigma_iy_i-bt},
\end{align*}
where $a_i$ and $b$ satisfy the following system of equations:
\begin{align*}
    &\sum_{i=1}^d a_i \Big(r - \frac{\sig^2_i}{2}\Big) + \frac 12\sum_{i=1}^d a_i^2\sig_i^2 + \sum_{i\ne j} \frac 12 \rho_{ij}\sig_i\sig_ja_ia_j-r-b = 0,\\
    &r-\frac{\sig_i^2}{2}+a_i\sig_i^2 +\sum_{j\ne i}\rho_{ij}\sig_i\sig_ja_j = 0,\qquad i = 1,\dots,d.
\end{align*}
Then, $u$ satisfies the following heat equation
\begin{align*}
\begin{cases}
\frac{\pl u}{\pl t} = \frac{1}{2}\sum_{i=1}^d \frac{\pl^2 u}{\pl y_i^2} + \sum_{i\ne j} \frac 12 \rho_{ij}\frac{\pl^2 u}{\pl y_i\pl y_j},\qquad (t,y)\in[0,T]\times\mathbb{R}^d,\\
u(0) = \Psi(e^{\sig_1y_1},\dots,e^{\sig_dy_d})e^{-\sum a_i\sig_iy_i}.
\end{cases}
\end{align*}
with $u(0) = \Psi(e^{\sig_1y_1},\dots,e^{\sig_dy_d})e^{-\sum a_i\sig_iy_i}$.

\begin{figure}[t]
\centering
\includegraphics[width=0.99\textwidth]{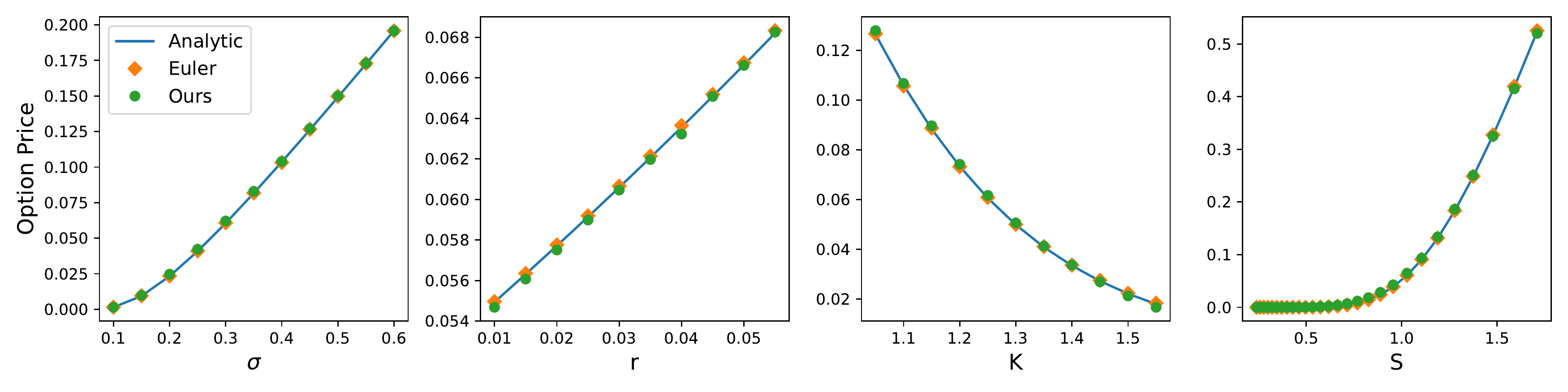}
\vspace{-2mm}
\caption{Ablation study on volatility $\sigma$, interest rate $r$, strike price $K$ and initial price. We fix a base setting with hyper-parameters $\sigma=0.03, r=0.3, K=1.25$, and run our algorithm on each setting with only one hyper-parameter deviated from the base setting. In addition, we plot the wave function under the base setting. We compare our solution at the execution time $T$ versus forward Euler method and the corresponding analytic ground truth. Our method is robust under all settings and achieves satisfactory performance.}
\label{fig:price}
\end{figure}

\subsection{Numerical Experiments}

In this section, we apply our algorithm to option pricing in Black-Scholes model across different settings using the numerical solution to the heat equation and the translation from the heat equation to the Black-Scholes equation from the previous subsection. We test the performance of our algorithm as well as show a calculation of the option price for higher dimensions.

In Figure \ref{fig:price} 
and 
Table \ref{tbl:black–scholes} we provide examples to test the performance of our algorithm. In both we include 1D examples for the price of a European call option, whose value at the expiry time $T$ is $V(T,s)=\Psi(s)=\max\{s-K,0\}$, where $K$ is a predetermined constant, called the strike price. We vary the volatility $\sigma$, strike price $K$, interest rate $r$, expiry time $T$, and initial price of the stock $S$. In Figure \ref{fig:price}, we compare against the forward Euler method and the ground truth {\it Black-Scholes formula}, which admits an analytical solution in 1D. Specifically, it is given by
\begin{align}
    V(t,s) = N(d_+)s-N(d_-)Ke^{-r(T-t)},
\end{align}
where $N$ is the cumulative distribution function of the standard normal distribution and 
$$
d_{\pm} = \frac{\ln{\frac{s}{K}}+(r\pm\frac{\sigma^2}{2})(T-t)}{\sigma\sqrt{T-t}}.
$$
Our method is robust under all settings and achieves satisfactory performance.
In Table \ref{tbl:black–scholes} we compare our accuracy against the forward Euler method, where we use the same analytical solution.  
Although the performance of our algorithm is inferior to that of Euler, it still achieves good accuracy and the margin can be treated as a price to pay for reducing the exponential scaling down to a polynomial one. Our approach is valuable for high dimensions, where other methods, such as the forward Euler method, suffer from the curse of dimensionality. Table \ref{tbl:black–scholes} also includes 2D examples with four different options: {\it Basket European call} and {\it put}, {\it rainbow max European call}, and {\it spread European put}, whose payoffs are listed in Table \ref{tbl:VarOpt:1}. Note that as we don't have an analytical solution for this case, the relative errors with respect to Euler solutions are reported instead.
In Figure \ref{fig:high_dim}, we provide a graph for the price of a basket European call option with up to five underlying stocks as a function of the strike price $K$. As expected the prices are convex with $K$.

\begin{table}[t]
\begin{minipage}[b]{0.5\linewidth}
\small
\centering
\begin{tabular}{ c|c }
Option Type & Payoff Function at expiry $\Psi(s)$\\
\hline
1D Call & $\max(s -K,0)$\\
\hline
Basket Call & $\max(\sum w_i s_i -K,0)$\\
Basket Put & $\max(K-\sum w_i s_i ,0)$\\
Rainbow Max Call & $\max(\max s_i -K,0)$\\
2D Spread Put & $\max(K-(s_1-s_2),0)$ \\
\end{tabular}
\vspace{6mm}
\caption{Payoff functions for our experiments. We consider basket call and put, Rainbow max call, and spread put options. 
}
\label{tbl:VarOpt:1}
\end{minipage}\hfill
\begin{minipage}[b]{0.46\linewidth}
\centering
\includegraphics[width=0.6\textwidth]{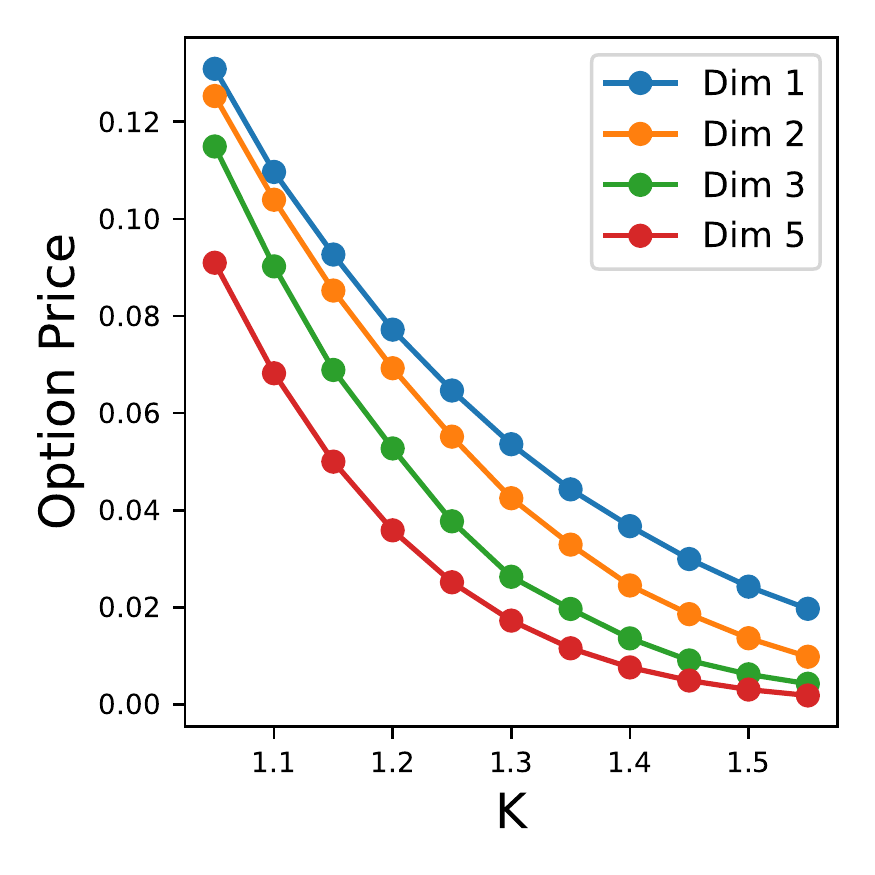}
\vspace{-3mm}
\captionof{figure}{Ablation study on dimensionality, following the base settings in Figure~\ref{fig:price}, $\sig = 0.3, r=0.03, K = 1.25$.} 
\label{fig:high_dim}
\end{minipage}
\end{table}

\begin{table*}[t]
\scriptsize
\centering
\caption{List of experiments for the application of our algorithm to Black–Scholes equation. The hyper-parameters for the experiments are listed. We compute the relative error of our method (\textbf{Ours}) at expiration time $T$ with analytical ground truth in the 1D case and Euler solution in the 2D case, respectively. For 1D, we also report the relative error of the forward Euler method with respect to the analytical ground truth (\textbf{Euler}) for comparison. Our algorithm achieves robust performance across various settings.}
\begin{tabular}{lc|c|c|cccc|cc}
\toprule
Problem & Initial \& Boundary Cond. & $d$ & $D$ & $T$ & $r$ & $K$ & $\sigma$ & Ours & Euler \\
\midrule
Black–Scholes 1D & 1D CALL & 1 & - & 1 & 0.03 & 1.25 & 0.3 & 0.011781 & 0.002494 \\
\midrule
Black–Scholes 1D & 1D CALL & 1 & - & 1 & 0.03 & 1.25 & 0.1 & 0.032792 & 0.000930 \\
Black–Scholes 1D & 1D CALL & 1 & - & 1 & 0.03 & 1.25 & 0.2 & 0.017272 & 0.002389 \\
Black–Scholes 1D & 1D CALL & 1 & - & 1 & 0.03 & 1.25 & 0.4 & 0.011560 & 0.001392 \\
\midrule
Black–Scholes 1D & 1D CALL & 1 & - & 1 & 0.03 & 1.05 & 0.3 & 0.086764 & 0.001750 \\
Black–Scholes 1D & 1D CALL & 1 & - & 1 & 0.03 & 1.15 & 0.3 & 0.013365 & 0.002521 \\
Black–Scholes 1D & 1D CALL & 1 & - & 1 & 0.03 & 1.35 & 0.3 & 0.012260 & 0.000328 \\
Black–Scholes 1D & 1D CALL & 1 & - & 1 & 0.03 & 1.45 & 0.3 & 0.012626 & 0.000845 \\
\midrule
Black–Scholes 1D & 1D CALL & 1 & - & 1 & 0.01 & 1.25 & 0.3 & 0.010436 & 0.002538 \\
Black–Scholes 1D & 1D CALL & 1 & - & 1 & 0.02 & 1.25 & 0.3 & 0.010599 & 0.002515 \\
Black–Scholes 1D & 1D CALL & 1 & - & 1 & 0.04 & 1.25 & 0.3 & 0.014192 & 0.002477 \\
Black–Scholes 1D & 1D CALL & 1 & - & 1 & 0.05 & 1.25 & 0.3 & 0.017423 & 0.002463 \\
\midrule
Black–Scholes 1D & 1D CALL & 1 & - & 0.5 & 0.03 & 1.25 & 0.3 & 0.022481 & 0.018385 \\
Black–Scholes 1D & 1D CALL & 1 & - & 1.5 & 0.03 & 1.25 & 0.3 & 0.012049 & 0.014881 \\
\midrule
Black–Scholes 2D & 2D BASKET CALL & 2 & 0.1 & 1 & 0.03 & 1.25 & 0.3 & 0.053477 & - \\
Black–Scholes 2D & 2D BASKET PUT & 2 & 0.1 & 1 & 0.03 & 1.25 & 0.3 & 0.043926 & - \\
Black–Scholes 2D & 2D RAINBOW MAX CALL & 2 & 0.1 & 1 & 0.03 & 1.25 & 0.3 & 0.057949 & - \\
Black–Scholes 2D & 2D SPREAD PUT & 2 & 0.1 & 1 & 0.03 & 1.25 & 0.3 & 0.031574 & - \\
\toprule
\end{tabular}
\label{tbl:black–scholes}
\end{table*}

Note that the Black-Scholes PDE lives on the positive orthant of $\mathbb{R}^d$ while its translation to the heat equation lives on $\mathbb{R}^d$. Thereby,
all numerical algorithms, including ours, require artificial truncation of the domain. We choose the hypercube domain to be $[s_l,s_u]^d=[K e^{-3\sigma_i\sqrt{T}},K e^{3\sigma_i\sqrt{T}}]^d$. 
This choice implies that $s_l$ is small (close to $0+)$ and $s_u$ is large (close $+\infty$).  On the faces of the hypercube, we use the time-discounted payoff functions, as they are reasonably accurate approximations of boundary values of the options considered. 
Given the number of qubits $n$ and the hypercube input domain for the heat equation $[L_l,L_u]^d$, which is approximately $[-5,5]^d$, the mesh size of each axis is $(L_u-L_l) / (2^{n/d}+1)$.

\section{Conclusions}
\label{sec:conclusions}
In summary, we introduced a generalization of McLachlan's variational principle applicable to generic time-dependent PDEs as well as a quantum-inspired training algorithm based on neural-network quantum states which can be used to perform approximate time evolution in high dimensions, overcoming the curse-of-dimensionality. Although we focused on a mesh-based formulation in which the quantum state vector is represented by $n$ qubits, it is clear that the mesh is not mandated by the formulation and it would be very interesting to pursue meshless variants based on continuous-variable neural-network quantum states including normalizing flows \cite{stokes2022numerical} and to address non-trivial boundary conditions. There exist a number of directions in which the results of this paper can be potentially improved. Since we only considered a first-order Euler approximation of the ODE \eqref{e:general} it would be natural to incorporate high-order time stepping schemes (e.g., Runge-Kutte methods). As an alternative, it would be interesting to pursue a direct solution of the discrete-time dynamical system \eqref{e:projected} which has proven successful in both the VMC \cite{gutierrez2022real} and VQA \cite{barison2021efficient} literature.
\section{Acknowledgements}

Authors gratefully acknowledge support from NSF under grants DMS-2038030 and DMS-2006305. This research was supported in part through computational resources and services provided by the Advanced Research Computing (ARC) at the University of Michigan.

\begin{appendix}

\section{Derivation of evolution equations}\label{app:evolution}
Assume that the time evolution map and the variational trial function admit Taylor expansions of the form
\begin{align}
    \Phi^{t+\delta t}_t\big(u\big) & = u + \mathcal{F}\big(t,u\big) \delta t + \mathcal{O}(\delta t^2) \enspace , \\
    u_{\theta + \delta\theta}
    & = u_\theta + \sum_{i=1}^p\frac{\partial u_\theta}{\partial \theta_i} \delta \theta_i
    + \mathcal{O}(\delta \theta^2) \enspace .
\end{align}
Then
\begin{align}
\left\Vert \Phi_{t}^{t+\delta t}(u_{\theta}) - u_{\theta + \delta\theta} \right\Vert_2^2
& = \sum_{i,j=1}^p \left\langle \frac{\partial u_\theta}{\partial \theta_i}  \middle| \frac{\partial u_\theta}{\partial \theta_j} \right\rangle \delta \theta_i \delta \theta_j  - \sum_{i=1}^p\left[ \left\langle \mathcal{F}(t,u_\theta) \middle| \frac{\partial u_\theta}{\partial \theta_i} \right\rangle + \left\langle \frac{\partial u_\theta}{\partial \theta_i} \middle| \mathcal{F}(t,u_\theta) \right\rangle \right] \delta t \, \delta \theta_i + \cdots  \notag \\
& = \sum_{i,j=1}^p \frac{1}{2}\left[\left\langle \frac{\partial u_\theta}{\partial \theta_i}  \middle| \frac{\partial u_\theta}{\partial \theta_j} \right\rangle + \left\langle \frac{\partial u_\theta}{\partial \theta_j}  \middle| \frac{\partial u_\theta}{\partial \theta_i} \right\rangle\right] \delta \theta_i \delta \theta_j  \notag \\
& \qquad \qquad \qquad \qquad - \sum_{i=1}^p\left[ \left\langle \mathcal{F}(t,u_\theta) \middle| \frac{\partial u_\theta}{\partial \theta_i} \right\rangle + \left\langle \frac{\partial u_\theta}{\partial \theta_i} \middle| \mathcal{F}(t,u_\theta) \right\rangle \right] \delta t \, \delta \theta_i + \cdots  \notag \\
& = \sum_{i,j=1}^p M_{ij}(\theta) \delta \theta_i \delta \theta_j  - 2\delta t \sum_{i=1}^p\delta \theta_i V_i(t,\theta) + \cdots
\end{align}
where in the last line we used the conjugate-symmetry of $\langle\cdot|\cdot \rangle$ and we have neglected $\delta\theta$-independent terms and terms higher than quadratic order in the multi-variable Taylor expansion in $\delta\theta$ and $\delta t$. The first-order optimality condition $0=\frac{\partial}{\partial \delta \theta_i}\left\Vert \Phi_{t}^{t+\delta t}(u_{\theta}) - u_{\theta + \delta\theta} \right\Vert_2^2$, gives, at lowest order in $\delta\theta$ and $\delta t$,
\begin{equation}
    0 = 2\sum_{j=1}^p M_{ij}(\theta)\delta \theta_j - 2 V_i(t,\theta) \delta t + \cdots
\end{equation}
and thus taking the limit $\delta t \longrightarrow 0$ gives the result.

\section{Matrix representation of $\mathcal{L}$}
\label{app:centraldiff}
In $d$ spatial dimensions and multi-index $\mathbf{i} = \{i_1,i_2,\dots,i_d\}$, let $\mathbf{i}\pm\mathbf{e}_k = \{i_1,i_2,\dots,i_k\pm 1,\dots,i_d\}$ and $\mathbf{i}\pm\mathbf{e}_k\pm\mathbf{e}_{k'} = \{i_1,i_2,\dots,i_k\pm 1,\dots,i_{k'}\pm 1,\dots,i_d\}$. Notice we do not allow $+1$ if $i_k = \frac nd$ or $-1$ if $i_k = 1$. Then the elements of the matrix is given by:
\begin{align}
[\widehat{\mathcal{L}}]_{\mathbf{i},\mathbf{j}} =
\begin{cases}
    \frac{-d}{\Delta^2}, \quad\mathbf{j}=\mathbf{i},\\   \frac{1}{2\Delta^2},\quad\mathbf{j}=\mathbf{i}\pm\mathbf{e}_k,\\
    \frac{D}{4\Delta^2},\quad\mathbf{j}=\mathbf{i}\pm\mathbf{e}_k\pm\mathbf{e}_{k'},\\
    \frac{-D}{4\Delta^2},\quad\mathbf{j}=\mathbf{i}\pm\mathbf{e}_k\mp\mathbf{e}_{k'},\\
    0, \quad \text{otherwise}.
\end{cases}
\end{align}
\end{appendix}

\bibliographystyle{plain}
\bibliography{refs}
\end{document}